\documentclass[11pt]{article}
\usepackage{amssymb}
\usepackage{amsmath}
\usepackage{fullpage}
\usepackage{latexsym,amsmath}
\usepackage{amsmath}

\usepackage{amsfonts}
\usepackage[english]{babel}
\usepackage[dvips]{graphicx}
\usepackage{graphics}
\usepackage{makeidx}
\usepackage[dvips]{color}

\usepackage{amsmath,amssymb}
\def \E {\mathbb{E}}
\def \Et {\mathbb{E}^{\mathcal{F}_t}}
\def \P {\mathcal{P}}
\def \calp {\mathcal{P}}
\def \e {\varepsilon}

\def\F {\mathcal{F}}
\def \R {\mathbb{R}}
\newtheorem{theorem}{Theorem}[section]
\newtheorem{lemma}[theorem]{Lemma}
\newtheorem{proposition}[theorem]{Proposition}

\numberwithin{equation}{section}
\def\Dim{\noindent\hbox{{\bf Proof.}$\;\; $}}          
\def\finedim{{\hfill\hbox{\enspace${ \square}$}} \smallskip}    
\def\sqr#1#2{{\vcenter{\vbox{\hrule height .#2pt
     \hbox{\vrule width .#2pt height#1pt \kern#1pt \vrule
     width .#2pt} \hrule height .#2pt}}}}
\def\square{\mathchoice\sqr54\sqr54\sqr{4.1}3\sqr{3.5}3}

\begin{document}
\title{Stochastic Maximum Principle for SPDEs with noise and control on the boundary}
\author{Giuseppina Guatteri, \\
Dipartimento di Matematica, \\
Politecnico di Milano, \\
Piazza Leonardo da Vinci 32, \\
20133 Milano, \\
Italia. \\
 {\tt e-mail: giuseppina.guatteri@polimi.it}}
\date{}
\maketitle
{\abstract
In this paper we prove necessary conditions for optimality of a stochastic control problem for a class of stochastic partial
 differential equations that is controlled through the boundary. This kind of problems can be interpreted as a stochastic control problem for an evolution system in a Hilbert space. The regularity of the solution of the adjoint equation, that is a backward stochastic equation in infinite dimension, plays a crucial role in the formulation of the maximum principle.  } 
 
\section{Introduction}
The maximum principle for stochastic control problems in infinite dimensions has been treated by Bensoussan  in \cite{Ben} using variational method. Then Hu and Peng in \cite{HuPeng1} studied general evolution controlled equations where the operator is a generator of a $C_0$- semigroup, see also the bibliography therein for the results in finite dimension and \cite{HuPeng2} for a stochastic equation of functional type.
There have been then several exstensions,  mainly dealing with finite dimensional systems and the few results regarding  SPDEs always consider diffused noise and control, for a comprehensive bibliography see \cite{YongZhou}.
The aim of our paper is to deal with control problems where the control and noise act on the boundary, a situation that seem interesting from the point of view of  applications.

As an example,  let us consider the following Cauchy problem for the stochastic heat equation
\begin{equation}\label{sistemaNeu_intro}\left\{\begin{array}{l}
\frac{\partial}{\partial t} v(t,x) = \frac{\partial^2}{\partial x^2} v(t,x) +
 f(v(t,x)) +
 g(v(t,x))  \frac{\partial^2}{\partial t\partial x} {\mathcal{W}}(t,x), \qquad  t \in [0,T], \ x \in [0,1]  \\ \\
\frac{\partial}{\partial x}v(t,0)= u^1(t) + \dot{W}^1(t), \qquad \frac{\partial}{\partial x}v(t,1)= u^2(t) + \dot{W}^2(t), \,\\ \\
v(0,x)= u^0(x),
\end{array}\right.
\end{equation}
where $\frac{\partial^2}{\partial t\partial x}{\mathcal{W}}(t,x)$ is a
space-time white noise, $\{ W^i_t, t \geq 0 \}$, $i=1,2$ are independent standard real Wiener processes, the unknown $v(t,x,\omega)$, representing the state of the system, is a real-valued process, the controls are two predictable real-valued processes $u^i(t,x,\omega)$, $i=1,2$ acting at $0$ and $1$, and $u_0$ is a function defined on $[0,1]$.

Equation \eqref{sistemaNeu_intro} can be rewritten as an evolution equation in an Hilbert space but to deal with the boundary terms one has to introduce unbounded terms in the equation.
Indeed if one sets $H=L^2(0,1)$ and $A$ the realization of the Laplace operator in $L^2$ with Neumann conditions one can write equation as 
\begin{equation}\label{forward_Neu_intro}
\left\{\begin{array}{ll}  dX_t = (A X_t + F(X_t))\, dt + (\lambda -A)D u_t \,dt + (\lambda -A)D_1\,d\tilde{W}_t + G(X_t)\, dW_t & t \in [0,T]\\
X_0 =x \end{array}\right.
\end{equation}
where $\lambda >0$ belongs to the resolvent set of the operator $A$, that is generator of an analytic semigroup and $D$ and $D_1$ are maps that transform the boundary terms $u$ and $\tilde {W}$
into elements that belong to the domain of fractional power of $\lambda -A$. Hence the operators $(\lambda -A)D $ and $(\lambda -A)D_1$ are the unbounded terms mentioned before and  are  regular enough to guarantee the existence of a {\em mild} solution  of \eqref{forward_Neu_intro} in the space $H$. 

A key point consists in proving that the solution to the adjoint equation is more regular and takes values in the domain of $[(\lambda -A)D]^*$ so that we can formulate the maximum principle. 

Since we do not assume in general the convexity of the control space we have to assume more regularity on the coefficients indeed  the rate of convergence of the first order approximations is $\e^{\alpha}$ with $ \frac{1}{2} <\alpha <1$.  
 The lower bound $ \frac{1}{2}$ is imposed by the presence of a noisy boundary term, in this case we prove a maximum principle condition without introducing the second order approximation of the state unknown and the additional adjoint equation. 

When the control space is a convex set, clearly this problem does not occur since the first order approximation is of order $\e$.

The rest of the paper is organized as follows: in section 2 we provide notation and we state the problem in his abstract formulation specifying the hypotheses, in section 3 we study the adjoint equation that turns out to be a backward stochastic equation in the infinite dimensional space $H$, in section 4 we prove the maximum principle, while in the last section we provide two examples of application of our result.
Notice that, exploiting the recent results of \cite{FabGol}, we can deal with  a heat equation with noisy boundary conditions of Dirichlet type. The drawback is that we have to work in an $L^2$ space with weight and so we have to restrict the class of the cost functionals we can treat.
\section{Preliminaries and statement of the problem}
\subsection{Notation}
Given a Banach space $X$,
the norm of its elements $x$ will be denoted
by $|x |_X$, or even by $|x|$ when no confusion is possible. If $V$ is
another Banach space, $L(X,V)$ denotes the space of bounded linear
operators from $X$ to $V$, endowed with the usual operator norm.
Finally we say that a mapping $F : X \to V$ belongs to the class $\mathcal{G}^1(X;V)$
if it is continuous, G\^ateaux differentiable on $X$,
 and $\nabla F: X\to L(X,V)$ is strongly continuous.
The letters $\Xi$, $H$, $K$ and $U$ will always be used to denote Hilbert spaces.
The scalar product is denoted $\langle \cdot, \cdot \rangle$, equipped with a
subscript to specify the space, if necessary. All the Hilbert
spaces are assumed to be real and separable; $L_2(\Xi,H)$ is the space of Hilbert-Schmidt operators from $\Xi $ to $H$, respectively.

\vspace{5pt}
 Given an arbitrary but fixed time horizon  $T$, we consider all stochastic processes
 as defined on subsets of the time interval $[0,T]$.
 Let $Q \in L(K)$ be a symmetric non-negative operator, not necessarily trace class and $\tilde{W} = (\tilde{W}_t)_{t \in [0,T]}$  be a $Q$-Wiener process with
values in $K$, defined on a
complete
probability space $(\Omega, \mathcal{F}, \mathbb{P})$ and $t{W} = ({W}_t)_{t \in [0,T]}$ be a cylindrical Wiener process with values in $\Xi$, defined on the same probability space and independent of $\tilde{W}$.
By $\{\mathcal{F}_t, \ t \in [0,T] \}$ we will denote the natural filtration
of $(\tilde{W},W)$, augmented with the family $\mathcal{N}$ of
$\mathbb{P}$- null sets of $\mathcal{F}$, see for instance \cite{DPZ1} for its definition. Obviously, the filtration
$(\mathcal{F}_t)$ satisfies the usual conditions of right-continuity and completeness. All the concepts
of measurability for stochastic processes will refer to this filtration.
By $\mathcal{P}$ we denote the predictable $\sigma$-algebra on
$\Omega \times [0,T]$ and by $\mathcal{B}(\Lambda)$ the Borel
$\sigma$-algebra of any topological space $\Lambda$.

Next we define two classes of stochastic processes with values in
a Hilbert space $V$.
\begin{itemize}
\item $L^2_\calp (\Omega\times [0,T];V)$ denotes the space of
equivalence classes of processes $Y \in L^2 (\Omega\times
[0,T];V)$ admitting a predictable version. It is endowed with the norm
\[ |Y|= \Big(\E \int_0^T |Y_s|^2 \, ds\Big)^{1/2}. \]
 \item $C_\calp([t,T];L^p(\Omega;S))$, $p\in [1,+\infty]$, $t \in [0,T]$, denotes the space of $S$-valued processes $Y$ such that
 $Y : [t,T] \to L^p(\Omega,S)$ is continuous
 and $Y$ has a predictable modification, endowed with the norm:
\begin{equation*}
|Y|^p_{C_\calp([t,T];L^p(\Omega;S))}=\sup_{s \in [t,T]}\E
|Y_s|^p_S
\end{equation*}
Elements of $C_\calp([t,T];L^p(\Omega;S))$ are identified up to
modification.
\item For a given $p \geq 2$, $L^p_{\mathcal{P}}(\Omega;C([0,T];V))$
     denotes the space of
    predictable processes $Y$ with continuous paths in $V$, such
    that the norm
    \[  \|Y\|_p  = (\E \sup _{s \in [0,T]} |Y_s|^p)^{1/p}\]
    is finite. The elements of $L^p_{\mathcal{P}}(\Omega;C([0,T];V))$
    are identified up to indistinguishability.
\end{itemize}
Given an element $\Phi$ of  $L^2_\P (\Omega\times
[0,T];L_2(\Xi,V))$ or of $L^2_\P (\Omega\times
[0,T];L_2(K,V))$, the It\^o stochastic integrals $\int_0^t
\Phi(s) \,dW(s)$ and  $\int_0^t \Phi(s) \,d\tilde{W}(s)$, $t \in
[0,T]$, are  $V$-valued martingales belonging to
$L^2_{\mathcal{P}}(\Omega;C([0,T];V))$. The previous definitions
have obvious extensions to processes defined on subintervals of
$[0,T]$ or defined on the entire positive real line $\R^+$.
\subsection{The optimal control problem and the state equation}
Let $H$ be a separable real Hilbert space, and $U$ a separable Hilbert, called the space of controls.
Let $U_{ad}$ a non empty set of $U$.
We set the space $L^2_\P (\Omega\times [0,T];U_{ad})$ the space of admissible controls, and we denote it by $\mathcal{U}$.
 
\noindent
We make the following, assumptions that we denote by ${\bf(A)}$:
\begin{enumerate}
\item[{\bf (A.1)}] $A: D(A)\subset H \to H$ is a 
linear, unbounded operator that generate a $C_0$-semigroup that is also analytic, $\{
e^{tA} \}_{t \geq 0}$  
such that $|e^{tA}|_{L(H,H)} \le M e ^{ \omega t}$, $ t \geq 0$
 for some $M >0$ and $\omega \in \R$. This means in particular that every $\lambda > \omega $ belongs to the resolvent set of $A$. \\
\item[{\bf(A.2)}] $F: \mathbb{R}^+ \times H \to H$,
 $G:\mathbb{R}^+ \times H \to L(\Xi,H)$, are measurable functions
such that for $h=
F,G$, $t \to h(t,x,y)$ is continuous for every fixed $x \in H, y \in K$.
Furthermore, there are constants $L$, $\Delta$ and  $\gamma \in
[0,1/2[$ such that:
\begin{align}
&| F(t,x)-
F(t,u)| _K \leq L |x-u|_H \nonumber\\
&|e^{sA}[G(t,x)- G(t,u)]|_{L_2(\Xi,H)}\leq \frac{L}{(1\wedge
s)^\gamma}|x-u|_H ,\nonumber \\
&|F(t,0)|_K \leq  \Delta , \nonumber \\
\nonumber
& |e^{sA}G(t,x)|_{L_2(\Xi,H)}  \leq \frac{\Delta }{(1\wedge
s)^\gamma}( 1+|x|_H),\end{align}
 for every $x,u
\in H$ and $s,t \in \mathbb{R}^+$.
\item[{\bf (A.3)}] $F (t, \cdot) \in
\mathcal{G}^{1}(H;H)$;
 for every $s >0$, \
 $e^{sA}G(t,\cdot) \in \mathcal{G}^{1}(H;L_2(\Xi,H))$ and 
\begin{align}
&| F_x(t,x)-
F_x(t,u)| _K \leq L |x-u|_H \nonumber\\
&|e^{sA}[G_x(t,x)- G_x(t,u)]|_{L_2(\Xi,H)}\leq \frac{L}{(1\wedge
s)^\gamma}|x-u|_H ,\end{align} for every $x,y \in H$;
 \item[{\bf (A.4)}] There exists a continuous linear operator  $D: U \to D((\lambda -A )^\alpha)$ for some $\frac{1}{2} < \alpha <1$ and $\lambda > \omega$, see for instance \cite{Lunardi} or \cite{Pazy} for the definition of the fractional power of the operator $A$. 
\item[{\bf (A.5)}] There exists a linear operator $D_1: U \to H$ and there is a constant $0< \beta < \frac{1}{2} $ such that the following holds:
\[  |e^{tA}(\lambda-A)D_1\sqrt{Q}|_{L_2(K,H)} \leq \frac{C}{t^\beta} \] 
 \end{enumerate}
 for some $\lambda >0$.
\vspace{2em}
 
 \noindent 
We consider in the 
Hilbert space $H$ the stochastic differential equation for the unknown process $X_t, \ t \in [0,T]$:
\begin{equation}\label{forward}
\left\{\begin{array}{ll}  dX_t = (A X_t + F(t,X_t))\, dt + (\lambda -A)D u_t \,dt + (\lambda -A)D_1\,d\tilde{W}_t + G(t,X_t)\, dW_t & t \in [0,T]\\
X_0 =x \end{array}\right.
\end{equation} 
As usual, see also \cite{DPZ1}, we mean by {\em mild} solution to this equation a $(\F_t)$- predictable process $X_t, \ t \in [0,T]$  with continuous path in $H$ such that $\P$- a.s.
\begin{multline*}
X_t = e^{tA} x + \int_0^t e^{(t-s)A}F(s,X_s)\, ds + \int _0^t e^{(t-s)A}(\lambda -A)D u_s \,ds  \\+ \int _0^t e^{(t-s)A}(\lambda -A)D_1 \, d\tilde{W}_s   +\int_0^t e^{(t-s)A} G(s,X_s)\, dW_s, \quad t \in [0,T]\\
\end{multline*}
\begin{proposition}
Under the assumptions ${\bf (A)}$, for every $u \in \mathcal {U}$ there exists a unique process 
$X \in  C_\P([0,T];L^2(\Omega;H))$  mild solution of equation \eqref{forward}. 
\end{proposition} 
\Dim
The only point to check in order to perform the fixed point argument, see\cite{DPZ1}, theorem 7.6, or \cite{FuhTess}, proposition 3.2, is that processes $(\int _0^t e^{(t-s)A}(\lambda -A)D u_s \,ds)_{t \in [0,T]} $ and $(\int _0^t e^{(t-s)A}(\lambda -A)D_1 \, d\tilde{W}_s)_{t \in [0,T]}$
 belong to the space $C_\P([0,T];L^2(\Omega;H))$. 
 We have indeed:
 \[\sup_{t\in [0,T]}\E \Big| \int_0^t e^{(t-s)A} (\lambda I -A) D u_s \, ds\Big|^2  \leq C^2 T^{1-2\alpha} \E \int_0^T |u_s|^2 \, ds  < +\infty \]\ 
Moreover:
\[\sup_{t\in [0,T]}\E \Big|\int_0^t e^{(t-s)A} (\lambda I -A) D_1 \, d \tilde{W}_s \Big|^2 \leq \textit{} C^2\E \Big[\int_0^T \rho^{-2\beta} \, ds \Big]^2 < +\infty \] 
 where both $C$ is defined in ${\bf (A)}$. 
 \finedim
 
We associate to this state equation the following  cost functional:
\begin{equation}\label{funzionale}
J(x,u)= \E \int_0^T l(t,X_t,u_t) \, dt + \E h(X_T)
\end{equation}
where $l$ and $h$ verify {\bf(B)}:
\begin{enumerate}
\item [{\bf (B.1)}] $(i)$ $l: [0,T] \times H \times U \to \R$
is measurable and there exist a constant $L >0$ and a modulus of continuity $\bar{\omega}:[ 0,+\infty) \to [0,+\infty)$, such that:
\begin{equation}
|l(t,x,u)-l(t,x',u')| \leq (L|x-x'| + \bar{\omega}(\|u-u'\|_{U}))
\end{equation}
for all $t \in [0,T]$, $x,x' \in H$ and $u,u' \in U$.

$(ii)$ Moreover for all $ t\in [0,T]$ and all $u \in U$ $l(t,\cdot,u) \in \mathcal{G}^1(H;\R)$ such that
\begin{equation}
|l_x(t,x,u)-l_x(t,x',u)| \leq L|x-x'| 
\end{equation}
for all $t \in [0,T]$, $x,x' \in H$ and $u \in U$.

\item [{\bf (B.2)}] $(i)$ $h: H  \to \R$, is measurable and there exist a constant $L >0$ such that
\begin{equation}
|h(x)-h(x')| \leq L|x-x'| 
\end{equation}
for all $x,x' \in H$.

$(ii)$ Moreover $h \in \mathcal{G}^1(H;\R)$ and 
\begin{equation}
|h_x(x)-h_x(x')| \leq L|x-x'| 
\end{equation}
for all $x,x' \in H$.\end{enumerate}
The optimal control problem consists in minimizing $J$ over all $u \in \mathcal{U}$.

We will seek for necessary conditions fulfilled by an optimal couple, whenever it exists, $(\bar{X},\bar{u})\in C_\P([0,T];L^2(\Omega;H)) \times \mathcal{U}$ such that 
\begin{equation}\label{min_funz} \inf_{u \in \mathcal{U}} J(x,u)= \E \int_0^T l(t,\bar{X}_t,\bar{u}_t) \, dt + \E h(\bar{X}_T)  \end{equation}  
where $\bar{X}$ is the mild solution to:
\begin{equation}\label{forward_ottima}
\left\{\begin{array}{ll}  d\bar{X}_t = (A \bar{X}_t + F(t,\bar{X}_t))\, dt + (\lambda -A)D \bar{u}_t \,dt + (\lambda -A)D_1\,d\tilde{W}_t + G(t,\bar{X}_t)\, dW_t & t \in [0,T]\\
\bar{X}_0 =x \end{array}\right.
\end{equation}
\section{Regularity results for the Adjoint equation}
In this section we consider the following backward stochastic differential equation, the so-called {\em adjoint equation}:
\begin{equation}\label{backward}
\left\{\begin{array}{ll}  - d{Y}_t = (A^T {Y}_t + F_x(t,\bar{X}_t)^T Y_t)\, dt +   G_x(t,\bar{X}_t)^T Z_t\, dt  - l_x(t,\bar{X}_t,\bar{u}_t) \, dt  -{Z}_t \, dW_t  - \tilde{Z}_t \, d\tilde{W}_t & t \in [0,T]\\

Y_T = - h_x(\bar{X}_T) \end{array}\right.
\end{equation}
Thanks to hypotheses {\bf(A)} on the derivatives $F_x$ and $G_x$  and hypotheses {\bf (B)} on the derivatives $l_x$ and $h_x$ this equation is affine with uniformly bounded coefficients (in the linear part) and integrable forcing term and integrable final data. The generator $A$  is an unbounded operator but generates a $C_0$-semigroup, so existence and uniqueness for the solution in $L^2_\P(\Omega;C([0,T];H)) \times  L^2_\P ([0,T]\times \Omega; L_2(K\times \Xi,H) )$ to this equation is a well known result, see \cite{HuPeng}.
It remains to prove an extra regularity property for the $Y$ component. 
\begin{proposition}
Under assumptions {\bf (A)} and {\bf (B)} there exists a unique mild solution $(Y, Z)$  in $L^2_\P(\Omega;C([0,T];H)) \times  L^2_\P ([0,T]\times \Omega; L_2(K\times\Xi,H) )$. 
Moreover for every $t \in [0,T[$ and  $\mathbb{P}-a.s.$, $Y_t(\omega)$ belongs to the domain of $D^TA^T$.
\end{proposition}  
\Dim
The mild solution exists by \cite{HuPeng}, theorem 3.1 or \cite{GuaTes}, theorem 4.4 that is a couple $(Y, Z)= (Y, (\hat{Z},\tilde{Z}))$  in $L^2_\P(\Omega;C([0,T];H)) \times  L^2_\P ([0,T]\times \Omega; L_2(K\times\Xi,H) )$ such that:
\begin{multline*}
Y_t = -e^{(T-t)A^T} h_x(\bar{X}_T) + \int_t^T e^{(s-t)A^T}(F_x(s,\bar{X}_s)^T {Y}_s + G_x(s,\bar{X}_s)^T {Z}_s)\, ds  \\ - 
\int_t^T e^{(s-t)A^T}l_x(s,\bar{X}_s, \bar{u}_s) \,ds  -\int_t^T e^{(s-t)A^T} Z_s \, dW_s -\int_t^T e^{(s-t)A^T} \tilde{Z}_s \, d\tilde{W} _s, \quad t \in [0,T]\\
\end{multline*}
Let us now prove the regularity result. We have:
\begin{align}\label{rappresentazione} 
Y_t = \Et Y_t= & - e^{(T-t)A^T} \Et h_x(\bar{X}_T) +  \int_t^T e^{(s-t)A^T}\Et (F_x(s,\bar{X}_s)^T Y_s + G_x(s,\bar{X}_s)^T Z_s)\, ds \nonumber\\
&-
\int_t^T e^{(s-t)A^T}\Et l_x(s,\bar{X}_s, \bar{u}_s) \,ds  
\end{align}
We have to evaluate $|| D^TA^T Y_t||_{U} = \sup_{u \in U, ||u||=1} \langle D^TA^T Y_t, u\rangle_U= \sup_{u \in U, ||u||=1} \langle  Y_t,AD u\rangle_U$.
We have:
\begin{align*}
&| \langle \Et h_x(\bar{X}_T), A e^{(T-t)A}Du \rangle | \leq \Et|h_x(\bar{X}_T)| |A e^{(T-t)A} Du| \leq \frac{C}{(T-t)^{1-\alpha}}(1+ \Et|\bar{X}_T|);\\
&|\int_t^T \langle \Et(F_x(s,\bar{X}_s)^T Y_s + G_x(s,\bar{X}_s)^T Z_s),A e^{(s-t)A}D u \rangle \,ds| \leq  C\int_t^T\frac{\Et(|{Y}_s|+ |Z_s|)}{(s-t)^{1-\alpha}} \,ds \\ &\leq C \Big(\int_0^T  (\Et|Y_s|^2 +\Et |Z_s|^2) \, ds\Big)^{1/2} T^{2\alpha -1};  \\
&|\int_t^T \langle \Et l_x(s,\bar{X}_s,\bar{u}_s),A e^{(s-t)A}D u \rangle\,ds | \leq \int_t^T  \frac{C}{(s-t)^{1-\alpha}}(1+\Et |\bar{X}_s| +\Et|\bar{u}_s|) \, ds  \\&\leq C \Big(\int_0^T  (1+\Et|\bar{X}_s|^2 + \Et|\bar{u}_s|^2) \, ds\Big)^{1/2} T^{2\alpha -1}. 
\end{align*}
This implies that for some constant $C >0$ that depends on $T$ and the quantities defined in {\bf (A)} and {\bf (B)}:
\[ \E|| D^TA^T Y_t||_{U} \leq \frac{ C }{(T-t)^{1-\alpha}}( 1 + \sup_{t \in [0,T]}\E |X_t|^2 + \E\sup_{t \in [0,T]} |Y_t|^2 + \E \int_0^T |Z_t|^2 \, dt ) < +\infty.  \]

\finedim
\section{The Maximum Principle}
\subsection{Variation of the trajectory}
Let $(\bar{X},\bar{u})$ be an optimal couple of problem \eqref{forward} and \eqref{min_funz}.
Fix $v \in U_{ad}$ and $ \bar{t} \in [0,T]$ and for every $ 0< \varepsilon < T-\bar{t}$ define
\begin{equation}\label{spike}
u^\e (t)=\left\{ \begin{array}{ll} 
v & t \in E_{\e}:= [\bar{t}, \bar{t} + \e]; \\
\bar{u}(t) & t \notin E_{\e}
\end{array} \right.
\end{equation}
Let us consider the following equations:
\begin{equation}\label{forward_spike}
\left\{\begin{array}{ll}  dX^\e_t = (A X^\e_t + F(t,X^\e_t))\, dt + (\lambda -A)D u^\e_t \,dt + (\lambda -A)D_1\,d\tilde{W}_t + G(t,X^\e_t)\, dW_t & t \in [0,T]\\
X_0^\e =x \end{array}\right.
\end{equation}
and 
\begin{equation}\label{forward_Iordine}
\left\{\begin{array}{ll}  d\tilde{X}^\e_t = (A \tilde{X}^\e_t + F_x(t,\bar{X}_t)\tilde{X}^\e_t)\, dt + (\lambda -A)D (u^\e_t-\bar{u}_t) \,dt +  G_x(t,\bar{X}_t)\tilde{X}^\e_t\, dW_t & t \in [0,T]\\
\tilde{X}_0^\e =0 \end{array}\right.
\end{equation}
We have:
\begin{proposition}\label{fattor}
Under hypothesis {\bf (A)} for every $\e >0$ there exist a unique mild solution  $X^\e \in C_\P([0,T];L^2(\Omega;H)) $ of equation \eqref{forward_spike} and a unique solution $\tilde{X}^\e \in C_\P([0,T];L^2(\Omega;H))$ of equation \eqref{forward_Iordine}. Moreover  for all $ p \geq  1$:
\begin{equation}\label{stimap}\E \sup_{t \in [0,T]} |\tilde{X}^\e_t|^p  < +\infty.\end{equation}
\end{proposition}
\Dim
The existence and uniqueness of the solutions are guaranteed by theorem  7.6 of \cite{DPZ1}.
Let us now prove \eqref{stimap}.
We have 
\begin{equation*}
\tilde{X}^\e_t = \int_0^t e^{(t-s)A} F_x (s,\bar{X}_s) \tilde{X}^\e_s \, ds + 
\int_0^t e^{(t-s)A} G_x (s,\bar{X}_s) \tilde{X}^\e_s \, dW_s  + \int_0^t  e ^{(t-s)A} (\lambda I -A) D ( u^\e_s - \bar{u}_s) \, ds
\end{equation*}
so 
\begin{align*}
\E \sup_{0 \leq t \leq r} |\tilde{X}^\e _t|^p \leq  c(p) &\Big[  \E \sup_{0 \leq t \leq r} \Big|\int_0^t e^{(t-s)A} F_x (s,\bar{X}_s) \tilde{X}^\e_s \, ds \Big| ^p +    \E \sup_{0 \leq t \leq r} \Big| \int_0^t e^{(t-s)A} G_x (s,\bar{X}_s) \tilde{X}^\e_s \, dW_s\Big| ^p \\ & +
 \E \sup_{0 \leq t \leq r} \Big|\int_0^t  e ^{(t-s)A} (\lambda I -A) D ( u^\e_s - \bar{u}_s) \, ds \Big| ^p \Big ]\leq c(p)[I_1+ I_2 + I_3]
\end{align*}
For $I_1$ we have, thanks to hypotheses {\bf(A.1)} and {\bf{(A.2)}} there exists a constant $C$ depending on $T$, $p$ and the quantities in ${\bf (A)}$ such that:
\begin{align*}
I_1 \leq C\E \sup_{0 \leq t \leq r}\int_0^t\sup_{0 \leq \sigma \leq s} |\tilde{X}^\e|^p \, ds \leq 
C\int_0^r\E \sup_{0 \leq \sigma \leq s} |\tilde{X}^\e|^p \, ds
\end{align*}
Then, having that $|A e^{(t-s)A}D|_{L(U,H)} \leq \frac{C}{(t-s)^{1-\alpha}}$ for soma costant $C>0$, thanks to {\bf (A.4)} we have
\begin{align*}
I_3 \leq C^p \Big(\int_0^T \frac{1}{s ^{1-\alpha}} \, ds \Big)^p |v|^p_U= C^p T^{\alpha p}|v|^p_U
\end{align*}
Eventually to treat term $I_2$ we use the factorization method, see \cite{DPKZ}. Take $ p >2$ and $ \rho \in (0,1)$ such that 
$ \frac{1}{p} < \rho < \frac{1}{2}-\gamma$, and let $ c_\rho^{-1}= \int_s^t (t-\sigma)^{\rho -1} (\sigma -s)^{-\rho} \, d\sigma$.
Hence
\begin{align*}
&I_2=\E \sup_{0 \leq t \leq r}\Big| \int_0^t e^{(t-s)A} G_x (s,\bar{X}_s) \tilde{X}^\e_s \, dW_s\Big| ^p \\ &
=\E \sup_{0 \leq t \leq r}\Big| \int_0^t e^{(t-\sigma)A} (t-\sigma)^{\rho-1} \,d\sigma \Big[\int_0^\sigma 
e^{(\sigma -s)A} (\sigma -s)^{-\rho}G_x (s,\bar{X}_s) \tilde{X}^\e_s \, dW_s\Big]\Big| ^p \\ & \leq  \E\sup_{0 \leq t \leq r}
\Big( \int_0^t (t-\sigma)^{(\rho -1)q} \Big)^{p/q} \int_0^t \Big|\int_0^\sigma 
e^{(\sigma -s)A} (\sigma -s)^{-\rho}G_x (s,\bar{X}_s) \tilde{X}^\e_s \, dW_s\Big|^p  \, d\sigma
\\& \leq C \int_0^r  \Big(\int_0^\sigma 
(\sigma -s)^{-2(\rho+\gamma)} \, ds\Big) ^{p/2} \E \sup_{0 \leq s \leq \sigma} |\tilde{X}^\e_\sigma|^p\, d\sigma.
\end{align*}
for some constant $C>0$, depending on $T$, $p$ and the parameter defined in ${\bf(A)}$.

\noindent
So combining all these estimates together, we obtain that 
\begin{equation*}
\E \sup_{  0\leq t \leq r } |\tilde{X}^\e_t|^p \leq C \Big[
\int_0^r \E \sup_{0 \leq \sigma \leq s} |\tilde{X}^\e_\sigma|^p \, ds  + |v|_U^p T^{\alpha p}. \Big] 
\end{equation*}
Thus by Gronwall theorem we can conclude.
\finedim

We claim that:
\begin{proposition}\label{prop_conv_unif}
Under  hypothesis {\bf (A)} there is a constant $\delta >0 $ independent of $\e$ such that:
\begin{equation}\label{con_unif_stato}
\Delta^\e:= \E \sup_{0 \leq t \leq T}|X^{\e}_t - \bar{X}_t| \leq \delta \e^\alpha
\end{equation}
\end{proposition}
\Dim
We have that:
\begin{equation}
\left\{\begin{array}{ll}  d({X}^{\e}_t - \bar{X}_t) = [A({X}^{\e}_t -\bar{X}_t)  + 
F(t,{X}^{\e}_t) - F(t,\bar{X}_t)]\, dt + (\lambda -A)D (u^\e_t-\bar{u}_t) \,dt \\  
\qquad\qquad\qquad\quad +(G(t,{X}^{\e}_t)- G(t,\bar{X}_t)\, dW_t, \qquad\qquad\qquad\qquad t \in [0,T]\\
{X}_0^{\e} - \bar{X}_0 =0 &\end{array}\right.
\end{equation}
 Following proposition \ref{fattor}, we have for $p >2$ :
\begin{align*}
\E \sup_{ t \in [0,r]} | X^{\e}_t - \bar X_t |^p  \leq C & \Big[\int_0^r \E \sup_{ \sigma \in [0,s]} | X^{\e}_\sigma -\bar{ X}_\sigma |^p \, ds +
\sup_{ \bar{t} \in [0,T]} \Big(\int_{\bar{t}}^{\bar{t} +\e}  \frac{|v|_U}{(t-s)^{1-\alpha}} \, ds\Big)^p \Big] 
\end{align*}
where $C$ as usual is independent of $\e$ and $m$ and it is function of $T$, $p$, $|v|_U$ and the quantities introduced in ${\bf(A)}$.
Therefore by Gronwall lemma we can conclude. 
\finedim

We set
\[ \eta^{\e}_t : = X^{\e}_t - \bar{X}_t - \tilde{X}^{\e}_t,  \] 
and we end the section with the following result
\begin{proposition}\label{prop_resto}
Under  hypothesis {\bf (A)} there is a constant $\delta_1 >0 $ independent of $\e$ such that:
\begin{equation}\label{con_unif_stato_resto}
\Delta_1\e:= \E \sup_{0 \leq t \leq T}|\eta^{\e}_t| \leq \delta \e^{2\alpha}
\end{equation}
\end{proposition}
\Dim
We have:
\begin{equation}
\left\{\begin{array}{ll}  d\eta^{\e}_t =  d(X^{\e}_t - \bar{X}_t - \tilde {X}^{\e}_t)=  
A \eta^{\e}_t \, dt + 
F_x(t,\bar{X}_t)\eta^{\e}_t \, dt + 
G_x(t,\bar{X}_t)\eta^{\e}_t \, dW_t\\ + \Big[\displaystyle\int_0^1 (F_x (t, \bar{X}_t +\theta (X^{\e}_t - \bar{X}_t)) - F_x (t,\bar{X}_t)) (X^{\e}_t - \bar{X}_t )\, d \theta \Big] \, dt\\
 + \Big[\displaystyle\int_0^1 (G_x (t, \bar{X}_t +\theta (X^{\e}_t - \bar{X}_t)) - G_x (t,\bar{X}_t)) (X^{\e}_t - \bar{X}_t )\, d \theta \Big] \, dt, \\
X^{\e}_0 - \bar{X}_0=0 &\end{array}\right.\\
\end{equation}
 Notice that, thanks to {\bf(A.3)}, we have that for some $ \gamma \in [0,\frac{1}{2}[$
\[
|e^{sA}[G_x(t,x)- G_x(t,u)]|_{L_2(\Xi,H)}\leq \frac{L}{(1\wedge
s)^\gamma}|x-u|_H \]
for every $x,y
\in H$ and $s \in \mathbb{R}^+$.
Thus for every $p > 2$, as in proposition \ref{fattor}, we obtain that
\begin{align*}
\E \sup_{ t \in [0,r]} | \eta^{\e}_t|^p \leq C & \Big[ \E \sup_{ \sigma \in [0,T]} | X^{\e}_\sigma - \bar{X}_\sigma |^{2p}  +
\int_0^r \E \sup_{ \sigma \in [0,s]} | \eta^{\e}_\sigma |^p \, ds \Big]
\end{align*} 
and we conclude.
\finedim
\subsection{Main result}
Now we are able to state and prove the {\em Maximum Principle} for our optimal control problem.
\begin{theorem}\label{princmax}
Assume hypotheses {\bf (A)} and {\bf (B)}.
Let $(\bar{X},\bar{u})$ be a optimal pair of Problem \eqref{min_funz} and \eqref{forward}.
Then there exists a unique solution $(Y,Z) \in L^2_\P(\Omega;C([0,T];H)) \times  L^2_\P ([0,T]\times \Omega; L_2(\Xi,H) )$ 
of equation \eqref{backward} and
\begin{equation}\label{hamiltoniana}
H(t, \bar{X}_t,\bar{u}_t,Y_t) \geq H(t, \bar{X}_t,v,Y_t), \qquad \forall v \in U_{ad},  \ a.e. \ t \in [0,T], \ \mathbb{P}-a.s.
\end{equation}
where 
\begin{equation*}
H(t,x,v,p):=
   \langle D^*(\lambda -A)^*p,v\rangle_H -l(t,x,v), \qquad (t,x,v,p) \in [0,T] \times H \times U_{ad} \times D(D^*(\lambda -A)^*), \ \lambda > \omega
 \end{equation*}
\end{theorem}
\Dim
Since $(\bar{X},\bar{u})$ is optimal for every $\e>0$ and $x\in H$ we have:
\[  0 \leq J(x,u_\e) - J(x,\bar{u}) = \E \int_0^T (l(t,X^\e_t,u^\e_t) - l(t,\bar{X}_t, \bar{u}_t))\,dt + \E (h(X^\e_T) - h(\bar{X}_T))= I_1+I_2.\]
 Let us considet $I_1$, adding and substractiong we get:
\begin{align*}
&\E \int_0^T (l(t,X^\e_t,u^\e_t) - l(t,\bar{X}_t, \bar{u}_t))\,dt=
 \E \int_0^T (l(t, {X}^{\e}_t, {u}^\e_t) - l(t,\bar{X}_t, {u}^\e_t))\,dt \\  & +  
\E \int_0^T (l(t,\bar{X}_t, {u}^\e_t) - l(t,\bar{X}_t, \bar{u}_t))\,dt   = J_1+J_2
\end{align*} 
Let us concentrate on $J_1$, we have, thanks to propositions \ref{prop_conv_unif} and \ref{prop_resto},:
\begin{align*}
&\E \int_0^T(l(t, {X}^{\e}_t, {u}^\e_t) - l(t,\bar{X}_t, {u}^\e_t))\,dt =  \\&
\E \int_0^T \int_0^1 [l_x(t, \bar{X}_t + \theta (X^{\e}_t- \bar{X}_t),u^\e_t )- l_x(t, \bar{X} _t, u_t^\e)] (X^{\e}_t- \bar{X} _t)\, d\theta \, dt 
   \\& + \E \int_0^T [l_x(t, \bar{X}_t,u^\e_t )- l_x(t, \bar{X} _t, \bar{u}_t)] (X^{\e}_t- \bar{X} _t)\,  \, dt+ \E \int_0^T l_x(t, \bar{X} _t,\bar{u}_t)\eta^{\e}_t \, dt  \\ & +\E \int_0^T l_x(t, \bar{X} _t, \bar{u}_t)\tilde {X}^{\e}_t \, dt  \leq C \e^{2\alpha} 
  +\E \int_0^T l_x(t, \bar{X} _t, \bar{u}_t)\tilde {X}^{\e}_t \, dt
\end{align*}
Combining all these relations we end up with:
\begin{align*}
&\E \int_0^T (l(t,X^\e_t,u^\e_t) - l(t,\bar{X}_t, \bar{u}_t))\,dt \leq C  \e^{2\alpha} + E \int_0^T l_x(t, \bar{X}_t, \bar{u}_t)\tilde{X}^{\e}_t \, dt+
 \E \int_0^T (l(t,\bar{X}_t, {u}^\e_t) - l(t,\bar{X}_t, \bar{u}_t))\,dt.
\end{align*}
Similarly we get:
\begin{align*}
\E h(X^\e_T) - \E h(\bar{X}_T) \leq C \e^{2\alpha} + \E h_x(\bar{X}_T)\tilde{X}^\e_T
\end{align*}
and thus:
\begin{align*}
&\E \int_0^T (l(t,X^\e_t,u^\e_t) - l(t,\bar{X}_t, \bar{u}_t))\,dt +\E h(X^\e_T) - \E h(\bar{X}_T) \\ & \leq C \e^{2\alpha}  
  + E \int_0^T l_x(t, \bar{X}_t, \bar{u}_t)\tilde{X}^{\e}_t \, dt+ \E h_x(\bar{X}_T)\tilde{X}^\e_T+ 
 \E \int_0^T (l(t,\bar{X}_t, {u}^\e_t) - l(t,\bar{X}_t, \bar{u}_t))\,dt
\end{align*}
Now, computing $d \langle \tilde{X}^\e_t, Y_t \rangle$, solutions respectively of equations \eqref{forward_Iordine} and \eqref{backward} we obtain:
\begin{align*}
-\E\langle \tilde{X}^\e_T,h_x( \bar{X}_T)  \rangle =  \E \int_0^T l_x(t, \bar{X}_t, \bar{u}_t)\tilde{X}^{\e}_t \, dt+ \E \int_0^T \langle AD (u^\e_t- \bar{u}_t), Y_t \rangle \, dt
\end{align*}
Therefore:
\begin{align*}
&0 \leq J(x,u^\e) - J(x, \bar{u})= \E \int_0^T (l(t,X^\e_t,u^\e_t) - l(t,\bar{X}_t, \bar{u}_t))\,dt+ \E h(X^\e_T) - \E h(\bar{X}_T) \\ & \leq C \e^{2\alpha}  
  + \E \int_{\bar{t}}^{\bar{t}+\e}\langle  (\bar{u}_t-v), D^*(\lambda -A)^*Y_t \rangle \, dt+ 
 \E \int_{\bar{t}}^{\bar{t}+\e} (l(t,\bar{X}_t, v) - l(t,\bar{X}_t, \bar{u}_t))\,dt
\end{align*}
Then dividing by $\e$ and recalling that $\alpha > \frac{1}{2}$,  
using a localization procedure, see \cite{Ben1}, and exploiting the continuity of $l$
 we can conclude.
\finedim
\subsection{The case when $U_{ad}$ is a convex set}
In this paragraph we assume that the non empty subset $U_{ad}$ is a {\em convex set}.
The space $\mathcal{U}$ of admissible controls and the optimal control problem are the same.

Besides ${\bf (A.1)}-{\bf (A.4)}- {\bf(A.5)}$ we assume ${\bf(C)}$:
\begin{enumerate}
\item[{\bf(C.1)}] $F: \mathbb{R}^+ \times H \to H$,
 $G:\mathbb{R}^+ \times H \to L(\Xi,H)$, are measurable functions
such that for $h=
F,G$, $t \to h(t,x)$ is continuous for every fixed $x \in H$.
Moreover $F (t, \cdot) \in
\mathcal{G}^{1}(H;H)$;
 for every $s,t >0$, \
 $e^{sA}G(t,\cdot) \in \mathcal{G}^{1}(H;L_2(\Xi,H))$ and there is a constant $L >0$ and $\gamma \in [0, \frac{1}{2}[$ such that:
 \begin{align}
&| F_x(t,x)| _K \leq L \nonumber\\
&|e^{sA}G_x(t,x)|_{L_2(\Xi,H)}\leq \frac{L}{(1\wedge
s)^\gamma},\end{align} for every $x \in H$ and every $ s,t \in \R^+$;
 
\item[{\bf(C.2)}] For all $ t\in [0,T]$ and all $u \in U$, $l(t,\cdot,u) \in \mathcal{G}^1(H;\R)$ and for all $ t\in [0,T]$ and all $x \in H$, $l(t,x,\cdot) \in \mathcal{G}^1(U;\R)$ and there is a constant $\Delta >0$ such that:
\begin{equation}
|l_x(t,x,u)|+ |l_u(t,x,u)| \leq \Delta ( 1+ |x|_H + |u|_U) 
\end{equation}
for all $t \in [0,T]$, $x \in H$ and $u \in U$.
\item[{\bf(C.3)}] $h \in \mathcal{G}^1(H;\R)$ and there is a constant $\Delta >0$ such that:
\begin{equation}
|h_x(x)| \leq \Delta ( 1+ |x|_H ) 
\end{equation}
for all $x \in H$.
\end{enumerate}

\vspace{2em} 

We have the following variational inequality.
\begin{lemma}
The cost functional $J$ is Gateaux-differentiable and the following variational inequality holds:
\begin{equation}\label{cond_I_ordine}
\frac{d}{d \theta} J(u(\cdot) + \theta v(\cdot))|_{\theta =0} = \E h_x(\bar{X}_T) \tilde{X}_T + \E \int_0^T [l_x(s, \bar{X}_s, \bar{u}_s) \tilde{X}_s +  l_u(s, \bar{X}_s, \bar{u}_s ) v_s] \, ds \geq 0
\end{equation} 
where $ v \in L_\P^2((0,T);U)$ satisfies $\bar{u}(\cdot) + v(\cdot) \in \mathcal{U}$, and $\tilde{X}$ is the solution to the following linear equation:
\begin{equation}\label{prima_approx}
\tilde{X}_t= \int_0^t  e ^{A(t-s)} F_x (s, \bar{X}_s) \tilde{X}_s  \, ds + \int_0^t  e ^{A(t-s)} G_x (s, \bar{X}_s) \tilde{X}_s  \, dW_s  + \int_0^t (\lambda -A) D v_s \, ds
\end{equation}
\end{lemma}

\Dim
The proof is similar to \cite{HuPeng1}.

\finedim 

It is clear that the results of Proposition 3.1 still hold under these hypotheses, so we can assert the maximum principle:
\begin{theorem} \label{massimo_convesso}
Let $(\bar {u},\bar{X})$ be an optimal pair the problem \eqref{min_funz} and \eqref{forward}. Then there 
exists a unique pair $(Y,Z) \in L^2_\P(\Omega;C([0,T];H)) \times  L^2_\P ([0,T]\times \Omega; L_2(\Xi,H) )$ solution of equation \eqref{backward} such that: 
\begin{equation}\label{hamiltoniana_conv}
\langle H_u(t, \bar{X}_t,\bar{u}_t,Y_t), v -\bar{u}_t \rangle \leq 0, \qquad \forall v \in U_{ad},  \ a.e. \ t \in [0,T], \ \mathbb{P}-a.s.
\end{equation}
where 
\begin{equation*}
H(t,x,u,p):=
   \langle D^*(\lambda -A)^*p,u\rangle_H -l(t,x,u), \qquad (t,x,u,p) \in [0,T] \times H \times U_{ad} \times D(D^*(\lambda -A)^*), \ \lambda > \omega
\end{equation*}
\end{theorem}
\Dim
From \eqref{cond_I_ordine} we have that:
\begin{equation} \label{cond_I}\E \langle h_x(\bar{X}_T), \tilde{X}_T\rangle  + \E \int_0^T \langle l_x(s, \bar{X}_s, \bar{u}_s), \tilde{X}_s \rangle + \langle  l_u(s, \bar{X}_s, \bar{u}_s ) ,v_s \rangle \, ds \geq 0
\end{equation} 
for every $ v \in L_\P^2((0,T);U)$ satisfies $\bar{u}(\cdot) + v(\cdot) \in \mathcal{U}$
 and $\tilde{X}$ is the solution to \eqref{prima_approx}.
 Moreover evaluating $\int_0^T d \langle  Y_t, \tilde{X}_t \rangle \, dt$, we get that:
 \begin{equation}\label{rel_fond}
 -\E \langle h_x(\bar{X}_T), \tilde{X}_T\rangle - \E \int_0^T \langle l_x(s, \bar{X}_s, \bar{u}_s) , \tilde{X}_s \rangle
= \E \int_0^T \langle Y_s, (\lambda -A)D v_s \rangle \, ds
 \end{equation}
Thus combining \eqref{cond_I} and \eqref{rel_fond} we end up with:
\[ \E \int_0^T  \langle l_u(s, \bar{X}_s, \bar{u}_s ), v_s \rangle  \, ds 
\geq   \E \int_0^T \langle D^*(\lambda -A)^*Y_s,  v_s \rangle \, ds 
\]
Thus using a localization procedure, see \cite{Ben1}, we can conclude chosing $v_t= v- \bar{u}_t$, with $v \in U_{ad}$. 
 \finedim
\section{Applications}
We provide two examples to which our result apply.
\subsection{Example 1} 
Let us consider the following problem:
\begin{equation}\label{sistemaNeu}\left\{\begin{array}{l}
\frac{\partial}{\partial t} y(t,x) = \frac{\partial^2}{\partial x^2} y(t,x) +
 f(y(t,x)) +
 g(y(t,x))  \frac{\partial^2}{\partial t\partial x} {\mathcal{W}}(t,x), \qquad  t \in [0,T], \ x \in [0,1]  \\ \\
\frac{\partial}{\partial x}y(t,0)= u^1(t) + \dot{W}^1(t), \qquad \frac{\partial}{\partial x}y(t,1)= u^2(t) + \dot{W}^2(t), \,\\ \\
y(0,x)= u^0(x),
\end{array}\right.
\end{equation}
where $\frac{\partial^2}{\partial t\partial x}{\mathcal{W}}(t,x)$ is a
space-time white noise, $\{ W^i_t, t \geq 0 \}$, $i=1,2$ are independent standard real Wiener processes, the unknown $y(t,x,\omega)$, representing the state of the system, is a real-valued process, the controls are two predictable real-valued processes $u^i(t,x,\omega)$, $i=1,2$ acting at $0$ and $1$ and having values in $\{-1, 0, 1\}$;
$u _0 \in L^2(0,1)$.
We assume that $f$ and $g$ are $C^{1,1} (\R)$. 

Now we write \eqref{sistemaNeu} as an evolution equation in the space $H=L^2(0,1)$. This is done in \cite{DPZ3}, see also bibliography therein and \cite{DebFuhTes}.
We define the operator $A$ in $H$ by setting 
\[D(A)=\{  y \in H^2(0,1): \frac{\partial}{\partial x}y(0)= \frac{\partial}{\partial x}y(1)=0 \}, \qquad Ay(x) = \frac{\partial^2}{\partial x} y (x), \text{ for  } y \in D(A). \]
Moreover for every $\lambda >0$,
\[ D((\lambda -A)^\alpha) = H^{2\alpha}(0,1) , \qquad \text{ for } 0 < \alpha < \frac{3}{4}.\]
For every fixed $\lambda >0$ there exists $d^i \in H^{2\alpha}(0,1)$, see for instance \cite{DebFuhTes} that solves the following Neumann problems:
 \begin{equation}\left\{\begin{array}{l}
\frac{\partial^2 }{\partial x^2}d^i (x) = \lambda d^i(x), \qquad x \in [0,1], \ i=1,2  \\ \\
\frac{\partial}{\partial x}d^1(0)= 1, \qquad \frac{\partial}{\partial x}d^1(1)=0, \,\\ \\
\frac{\partial}{\partial x}d^2(0)= 0, \qquad \frac{\partial}{\partial x}d^2(1)=1.\,\\ \\
\end{array}\right.
\end{equation}
Thus we set $U=K=\R ^2$ and $U_{ad} = \{-1, 0, 1 \} \times  \{-1, 0, 1 \}$, the covariance matrix $Q=I$ and $D=D^1:\R ^2\to D((\lambda -A)^\alpha)$, such that $D u(t)(x) = d^1(x) u^1(t) + d^2(x) u^2(t)$ and  
$D \tilde{W}(t)(x) = d^1(x) W^1_t + d^2(x) W^2_t$.
We set $\Xi=L^2(0,1)$ and $\frac{\partial}{ \partial x}{\mathcal{W}}(t,\cdot) = W(t)$ is a cylindrical Wiener process in
$ \Xi= L^2([0,1])$, see for instance \cite{DPZ1}. \\
We finally set $X_t = y(t,\cdot)$ and $F(\xi)(\cdot)= f(\xi(\cdot))$ and $G(\xi)(\cdot)= g(\xi(\cdot))$ for all $ \xi \in H$, then system \eqref{sistemaNeu}
can be written as
\begin{equation}\label{forward_Neu}
\left\{\begin{array}{ll}  dX_t = (A X_t + F(X_t))\, dt + (\lambda -A)D u_t \,dt + (\lambda -A)D_1\,d\tilde{W}_t + G(X_t)\, dW_t & t \in [0,T]\\
X_0 =u_0,  \end{array}\right.
\end{equation}
It is then easy to show that all hypotheses {\bf (A)} are fulfilled, note that we can chose $\alpha > \frac{1}{2}$.
Let us introduce the following finite horizon cost
\[ J(x,u^1(\cdot),u^2(\cdot)) = \E \int_0^T \int_0^1 \text{l}(y(s,x), u^1 (s), u^2 (s)) \, dx\, ds + \E \int_0^1 \text{h}(y(T,x))\, dx .\]
with $\text{l}(x,u^1,u^2) :\R\times \R \times \R \to \R$, continuous and derivable, with bounded and Lipschitz continuous derivatives in $x$ uniformly with respect to $u^1, u^2$ and Lipschitz continuous with respect to $(u^1,u^2)$ uniformly with respect to $x$. Moreover we assume $\text{h}\in C^{1,1}(\R)$.
We set: 
\[l(\xi,u^1,u^2)= \int_0^1  \text{l}(\xi(x), u^1, u^2) \, dx, \text{ for } \xi \in H, u^1, u^2 \in U\]
\[h(\xi)=\int_0^1 \text{h}(\xi(x))\, dx, \text{ for }  \xi \in H.   \]
Hence also hypothesis {\bf (B)} is fulfilled and theorem  \ref{princmax} holds.
\subsection{Example 2}
Now we consider a boundary control problem for a stochastic heat equation with Dirichlet condition perturbed by a stochastic process.
More precisely we have:
\begin{equation}\label{sistemaDir}\left\{\begin{array}{l}
\frac{\partial}{\partial t} v(t,x) = \frac{\partial^2}{\partial x^2} v(t,x) +
 f(v(t,x)), \qquad  t \in [0,T], \ x \in [0,+\infty)  \\ \\
v(t,0)= u(t) + \dot{W}(t),  \,\\ \\
v(0,x)= v^0(x),
\end{array}\right.
\end{equation}
where  $\{ W_t, t \geq 0 \}$, is a standard real Wiener process, the unknown $v(t,x,\omega)$, representing the state of the system, is a real-valued process, the control is a predictable real-valued process $u(t,x,\omega)$,  acting at $0$;
 $v^0 \in L^2(0,1)$.
We assume that $f$ is $C^{1} (\R)$ with bounded derivative, and clearly $g=0$. 
It is well known that it is not possible to rewrite the Cauchy problem \eqref{sistemaDir} as an evolution equation in the space $L^2(0,+\infty)$,  see \cite{DPZ3}. In \cite{Krylov} it is shown that the Dirichlet map takes values in the domain of $(-A)^\alpha$, for a certain $\alpha > \frac{1}{2}$, when it is considered in the Hilbert space $L^2(0,+\infty;(\rho^{\theta +1}\wedge 1) \, d\rho )$.
More precisely if we set $H = L^2(0,+\infty;(\rho^{\theta +1}\wedge 1) \, d\rho )$, the operator $A_0$ that is the realization of the Laplacian with Dirichlet conditions 
in $L^2(0,+\infty)$ extends to an operator $A$ that generates an analytic operatorin $H$.

For every fixed $\lambda >0$ there exists $d \in D((\lambda-A)^\alpha))$ for some $\alpha > \frac{1}{2}$:
\begin{equation}\left\{\begin{array}{l}
 \frac{\partial^2}{\partial x^2}d (x) = \lambda d(x), \qquad x \geq 0  \\ \\
d(0)= 1, 
\end{array}\right.
\end{equation}

Thus we set $U= U_{ad}=K=\R $, the covariance matrix $Q=1$ and $D=D^1:\R ^2\to D((\lambda-A)^\alpha))$, such that $D u(t)(x) = d(x) u(t) $ and  
$D \tilde{W}(t)(x) = d(x) W_t $.
We finally set $X_t = v(t,\cdot)$ and $F(\xi)(\cdot)= f(\xi(\cdot))$ for all $ \xi \in H$, then system \eqref{sistemaDir} can be written as an evolution equation in $H$. 

We assume the cost functional is of the following form:
\[l(\xi,u)= \int_0^{+\infty}  \frac{(x\wedge 1)^{1+\theta}}{(1+x^2)^{1/2+\varepsilon}}\text{l}(\xi(x), u) \, dx, \text{ for } \xi \in H, u \in U, \text{ for some }\varepsilon > 0\]
\[h(\xi)=\int_0^{+\infty}\frac{(x\wedge 1)^{1+\theta}}{(1+x^2)^{1/2+\varepsilon}} \text{h}(\xi(x))\, dx, \text{ for }  \xi \in H, \text{ for some }\varepsilon > 0\  \]
Where $\text {l} : \R\times \R \to \R$ is continuous and derivable  with continuous derivarives with sublinear growth and $h$ is derivable with the derivative with sublinear growth.
Thus theorem \ref{massimo_convesso} holds.

\end{document}